# *GeoGebra e situações que envolvem modelação numa abordagem STEAM*

José Manuel Dos Santos [1]

Astrigilda Pires Silveira [2]

Alexandre Emanuel da Silva Trocado[3]

## *Resumo:*

No sentido se implementar uma abordagem STEAM incluindo o uso da tecnologia, nomeadamente o uso do software de matemática interativa o GeoGebra, nas aulas de matemática, no espaço lusófono, foram concebidos os materiais que aqui se apresentam, a implementar numa primeira fase junto dos professores. Posteriormente, com as necessárias adaptações, estas tarefas serão aplicadas junto dos estudantes. As tarefas abordam situações de modelação, em problemas geométricos a duas e três dimensões, com o intuito de aplicar o software GeoGebra na sua análise de modo a ilustrar as suas capacidades. São utilizadas as diferentes janelas deste software, nomeadamente as janelas 2D e 3D, janela CAS, folha de cálculo e as janelas bidimensionais extra de modo a estudar planos de corte em sólidos e algumas superfícies. As tarefas são apresentadas para que qualquer utilizador, independente do grau de conhecimento que tenha do software, as possa acompanhar, suportando-se em roteiros com algumas indicações das ferramentas e comandos a usar. Pensadas para o ensino e aprendizagem da Matemática, a partir de uma abordagem STEAM, estas tarefas permitem conexões com outras Ciências e as Artes, e propiciam o desenvolvimento de projetos utilizando e consolidação conteúdos matemáticos relevantes. Estas tarefas integram as propostas de atividades dos participantes dos Cursos de Formação de Formadores em GeoGebra para os Países de Língua Oficial Portuguesa, que a partir de 2019 têm incidência na abordagem STEAM. Estes cursos são realizados com o alto patrocínio da Organização de Estados Ibero-Americanos para a Educação, a Ciência e a Cultura (OEI). Dado o interesse que as tarefas tem para os utilizadores do espaço ibérico, bem como a sua disseminação a nível global, os materiais que inicialmente foram elaborados em língua portuguesa serão adaptadas para os falantes de língua castelhana e inglesa.

## *Abstract:*

In order to implement a STEAM approach including the use of technology, namely the use of interactive mathematics software GeoGebra, in mathematics classes, in the lusophone space, the materials presented here were conceived, to be implemented in a first phase among teachers. Later, with the necessary adaptations, these tasks will be applied to the students. The tasks deal with modeling situations, in two- and three-dimensional geometric problems, in order to apply GeoGebra software in its analysis to illustrate its capabilities. The different windows of this software are used, namely the 2D and 3D windows, CAS window, spreadsheet and extra two dimensional windows in order to study cutting planes in solids and some surfaces. The tasks are presented so that any user, regardless of the degree of knowledge they have of the software, can follow them, being supported in scripts with some indications of the tools and commands to use. Designed for the teaching and learning of Mathematics, from a STEAM approach, these tasks allow connections with other Sciences and the Arts, and allow the development of projects using and consolidating relevant mathematical contents. These tasks are part of the proposals of activities of the participants of the Training Courses for Trainers in GeoGebra for Portuguese Speaking Countries, which from 2019 have an impact on the STEAM approach. These courses are carried out with the high sponsorship of the Organization of Ibero-American States for Education, Science and Culture (OEI). Given the interest that the tasks have for the users of the Iberian space, as well as their dissemination at a global level, the materials initially developed in Portuguese language will be adapted for Spanish and English speakers.

## *Abstract:*

Para implementar un enfoque STEAM que incluya el uso de la tecnología, a saber, el uso del software interactivo de matemáticas GeoGebra, en las clases de matemáticas, en el espacio lusófono, los materiales presentados aquí fueron concebidos, para ser implementados en una primera fase entre los maestros. . Posteriormente, con las adaptaciones necesarias, estas tareas serán aplicadas a los alumnos. Las tareas se refieren a situaciones de modelado, en problemas geométricos de dos y tres dimensiones, para aplicar el software GeoGebra en su análisis para ilustrar sus capacidades. Se utilizan las diferentes ventanas de este software, a saber, las ventanas 2D y 3D, la ventana CAS, la hoja de cálculo y las ventanas bidimensionales adicionales para estudiar los planos de corte en sólidos y algunas superficies. Las tareas se presentan para que cualquier usuario, independientemente del grado de conocimiento que tengan del software, pueda seguirlas, y sea compatible con scripts con algunas indicaciones de las herramientas y los comandos que debe utilizar. Diseñadas para la enseñanza y el aprendizaje de Matemáticas, desde un enfoque STEAM, estas tareas permiten conexiones con otras Ciencias y Artes, y permiten el desarrollo de proyectos utilizando y consolidando contenidos matemáticos relevantes. Estas tareas forman parte de las propuestas de actividades de los participantes de los Cursos de capacitación para capacitadores en GeoGebra para países de habla portuguesa, que desde 2019 tienen un impacto en el enfoque STEAM. Estos cursos se llevan a cabo con el alto patrocinio de la Organización de Estados Iberoamericanos para la Educación, la

[1] Instituto Geogebra Portugal, ESE-PP, AEDCM, santosdossantos@ese.ipp.pt .

[2] Instituto GeoGebra Cabo Verde, Universidade de Cabo Verde, astrigilda.silveira@docente.unicv.edu.cv .

[3] Instituto Geogebra Portugal, mail@alexandretrocado.com

Ciencia y la Cultura (OEI). Dado el interés que tienen las tareas para los usuarios del espacio ibérico, así como su difusión a nivel global, los materiales desarrollados inicialmente en idioma portugués se adaptarán para hablantes de español e inglés.

## *Tecnologia, Matemática e STEAM*

No ensino e aprendizagem da matemática contemporâneo, em todos os níveis de ensino, carece de integrar não só a tecnologia, mas também o estabelecimento de conexões com outras áreas do saber (NCTM, 2000), nomeadamente com as ciências em geral. Em relação a vantagem do uso da tecnologia, no ensino da matemática e nos seus efeitos no desenvolvimento profissional, nomeadamente no ensino básico e secundário, encontra-se bastante estudada. Vários autores (Arcavi & Hadas, 2000; Jones, 2000; Olive, 2000; Hohenwater, et al, 2008; Hollebrands, 2007; Heid & Blume, 2008; Haciomeroglu, et al, 2009; Caridade, 2012; Barrera-Mora & Reyes-Rodríguez, 2013; Tabaghi & Sinclair, 2013; Vitale, Swart, & Black, 2014, Tomaschko, Kocadere & Hohenwarter, 2018) são unânimes em afirmar que o ensino da Matemática, apoiado em actividades suportados pela tecnologia, favorecem o desenvolvimento de atitudes positivas que conduzirão a uma melhor aprendizagem e mais gosto por esta ciência. Mesmo no ensino superior,

> Modern science, technology, engineering and mathematics (STEM) education is facing fundamental challenges. Most of these challenges are global; they are not problems only for the developing countries. Addressing these challenges in a timely and efficient manner is of paramount importance for any national economy. Mathematics, as the language of nature and technology, is an important subject in the engineering studies. Despite the fact that its value is well understood, students' mathematical skills have deteriorated in recent decades in the western world. This reflects in students' slow progressing and high drop-out percentages in the technical sciences.
>
> The remedy to improve the situation is a pedagogical reform, which entails that learning contexts should be based on competencies, engineering students motivation should be added by making engineering mathematics more meaningfully contextualized, and modern IT-technology should be used in a pedagogically appropriate way so as to support learning. ( Pohjolainen, Myllykoski, Mercat, & Sosnovsky,2018)

O estudo anteriormente citado compara a situação da educação STEM nos EU, Europa, Rússia, Arménia e Geórgia, colocando como recomendação que os curricula tendam a privilegiar a matemática aplicada, modelando situações concretas e problemas, de modo a aumentar o nível de motivação dos estudantes engajando-os para a aprendizagem, à semelhança do que acontece na generalidade dos países europeus e no EUA. Já em relação a tecnologia, no mesmo estudo, são citados os programas de distribuição livre, o Sage e o GeoGebra, entre outros, sendo o GeoGebra usado em diferentes casos de estudo. Existe também a recomendação do incremento do uso de sistemas de gestão de aprendizagem (Learning management system) para apoiar a avaliação formativa e o trabalho autónomo dos estudantes, a semelhança do que já se faz em muitos dos países europeus (Pohjolainen, Myllykoski, Mercat, & Sosnovsky, 2018, p. 194).

No caso do ensino não superior o termo STEM refere-se a uma abordagem educativa onde a Matemática é colocada ao estudante em contexto, propondo-lhe a realização de projetos onde várias competências e conhecimentos são adquiridas com a utilização de múltiplas tecnologias. Esta abordagem está contemplada desde 1990 nos currículos americanos, tem vindo a ser amplamente usada em vários curriculas do Estados Unidos (Kelley & Knowles, 2016). Note-se também que a abordagem STEM, paulatinamente tem sido assumida em projectos realizados em vários países europeus. A abordagem STEM passou também a incluir as áreas das artes e do design, passando a ser designada pelo acrónimo STEAM (Bequette & Bequette, 2012).

Na visão STEM para 2026, subscrita em 2015, pelo grupo de peritos convocados pela Departamento de Educação dos EUA, em colaboração com os Institutos Americanos de Investigação (AIR), aponta para a promoção de experiências de aprendizagem e ensino de qualidade e culturalmente relevantes para todas as crianças e jovens. Nesta estratégia, os estudantes acedem e possuem um sentimento de pertença a uma comunidade STEM, comunidade que promove um percurso de aprendizagem ao longo da vida, que se desenvolve em ambientes formais e informais desde o ensino básico ao ensino superior, incluindo escolas, centros de ciência e outras instituições. A visão STEM 2026 integra seis componentes, fortemente inter-relacionadas: i) comunidades de prática engajadas e em rede; ii) atividades de aprendizagem acessíveis que convidam a recreação e ao assumir de riscos intencionais; iii) experiências educacionais que incluem abordagens interdisciplinares para resolver "desafios"; iv) a espaços de aprendizagem flexíveis e inclusivos suportados por tecnologias inovadoras (Boaler, 2016).; v) formas diversificadas de avaliação privilegiando a avaliação formativa e o feedback; vi) imagens e ambientes, societais e culturais, que promovem diversidade e oportunidades em STEM (Tanenbaum, 2016).

## *GeoGebra, a razão da opção*

De acordo com as orientações do National Council Teachers of Mathematics (2008), com recurso aos Ambientes Dinâmicos de Geometria Dinâmica (ADGD) os estudantes poderão aprender melhor a matemática. De entre os ADGD, destaca-se o GeoGebra, software matemático dinâmico abrangendo as quatro grandes áreas da Matemática – a Geometria, a Álgebra, o Cálculo e a Probabilidade e Estatística –, de cariz predominantemente construtivista, constitui um excelente recurso para o estudo da Matemática e de outras áreas (Química, Física, ...). Segundo Hohenwarter e Preiner (2007), a múltipla perceção dos objetos – por exemplo, cada expressão na zona algébrica possui uma representação na zona gráfica e vice-versa –, constitui a característica mais peculiar do GeoGebra, comparada com outros ambientes dinâmicos (Hohenwarter, 2013). A possibilidade de estudante ver, explorar, conjeturar, validar, compreender e comunicar os conceitos geométricos de uma forma interativa e atrativa encontra no GeoGebra, um recurso apropriado e moderno para o estudo da Matemática em diversos graus de ensino (Fernandes, 2018; Silveira, 2015). Saliente-se ainda, o papel mediador da Matemática com as ciências e tecnologia, dão ao software livre GeoGebra condições ímpares e que o investem como ferramenta de eleição para o ensino e aprendizagem das ciências em geral. Acresce ainda, a circunstância do GeoGebra ser um software gratuito é, para além de todas as valias já enunciadas, um aspeto muito importante para a realidade de muitos dos países de língua portuguesa.

O facto deste software dispor de diversas janelas ou vistas, que intercomunicam com o uso de álgebra computacional, permite-lhe um número infindável de aplicações a situações de modelação da ciência em geral, que usam a matemática como uma linguagem privilegiada na análise e descrição dos seus objetos de estudo, em todos os graus de ensino. Acresce ainda que a possibilidade de lidar com múltiplas janelas conferem ao GeoGebra características únicas que potenciam a imagética associada ao pensamento e raciocínio matemático, imagética está que tem um papel fundamental na procura de generalizações (Presmeg, 1997) e na potenciação do raciocínio abdutivo. A integração de diferentes janelas gráficas, 2D e 3D, que intercomunicam de modo dinâmico, conferem a estas janelas modelos das extensões "naturais" entre R2 e R3, colocando este aplicativo como um auxiliar inestimável na: visualização de relações estabelecidas em diversos problemas de modelação matemática; utilização de Calculo Algébrico e Simbólico (CAS); possibilidade de programação, nomeadamente na programação de autómatos; no desenvolvimento de atividades de descrição e inferência estatística; demonstração automática de teoremas em geometria; integração do software na generalidade dos sistemas de gestão de aprendizagens existentes.

Do ponto de vista da didática da disciplina de Matemática, o GeoGebra é utilizado por uma ampla comunidade internacional de professores, as suas aplicações em múltiplos contextos educativos tem sido objeto de inúmeros estudos focando-se na aprendizagem dos estudantes nos primeiros anos (da Silva & Fernandes, 2015, 2017; Furner & Marinas, 2012), no ensino básico (Silveira & Cabrita, 2013; Cadavez, Morais, Cadavez & Miranda, 2013), no ensino superior (Batalla, 2019, Breda & Dos Santos, 2013, 2015, 2016; Maia-Lima, Silva & Duarte, 2015 ; Sousa, Furtado & Horta, 2018; Turgut, 2019), e no ensino de outras ciências (Maximiano, Silveira, Belloti, Rossini, Bacilli, & Orlandi, 2012; Malgieri, Onorato & De Ambrosis, 2014; Marciuc & Miron, 2014; Solvang & Haglund, 2018 ). Considerando a investigação contemporânea, o software tem sido amplamente usado, na área da educação existe um grande número de estudos em didática, também múltiplos estudos relacionados com a formação de professores têm sido realizados (Zengin, 2019). Na investigação em Matemática recentemente o GeoGebra está a ser usado na área da geometria esférica (Breda & Dos Santos, 2018, 2019), mesmo na área da Física há estudos que tem usado o GeoGebra (Scheuber, 2018; Fatema & Singh, 2019).

Existe uma ampla comunidade global de utilizadores do GeoGebra, esta inclui estudantes, educadores e investigadores em várias áreas das ciências, e naturalmente investigadores de vários ramos da matemática. De facto entre os softwares de utilização livre criados para trabalhar a matemática e as suas aplicações em geral o

> "GeoGebra as a potentially defining moment in Kaput's grand vision of democratic access; others maintain that the software still presupposes technological resources which are not always available in financially weak regions; and still others feel that the open-source status of the software is unsustainable

> in terms of development, growth, and quality control. However, we contend that these latter objections are likely unfounded, particularly in light of other longstanding and successful open-source initiatives (e.g., Linux), and the fact that several large software companies are now developing products that involve open-source technology in order to meet the challenges of sustainability and relevancy within a rapidly changing marketplace." (Jarvis et al. , 2011, pp. 232-233).

Acresce à ideia de acesso a um sistema de aprendizagem de qualidade, a necessidade de integração da tecnologia nas diversas áreas de conhecimento. Neste ponto, a integração das diversas componentes do GeoGebra nos dispositivos moveis, nomeadamente nos smarthphones, transformam um objecto quotidiano numa calculadora gráfica, com capacidades CAS, e que pode funcionar em modo de exame, que de resto já é utilizado em exames nacionais da Áustria, em algumas regiões da Alemanha e na Noruega.

Refira-se ainda que, em Moçambique, no desenvolvimento de um projeto anterior da OEI, inserido no programa Iberciencia, a utilização de dispositivos móveis foi uma mais valia para levar as atividades com o GeoGebra as salas de aula, situação que também aconteceu em Cabo Verde nas experiências de ensino desenvolvidas no projeto que levou a instalação do Instituto GeoGebra na Universidade de Cabo Verde.

Apesar da dificuldade de infraestruturas físicas na rede das escolas africanas, é também um facto que em muitas escolas os estudantes possuem telefones moveis, com planos de dados incluídos, que lhes permite usar o GeoGebra em sala de aula, e note-se que estudantes de algumas escolas de Maputo se prontificaram a usar os seus dispositivos e dados moveis para desenvolver as tarefas que o professor levou para a sala de aula. Situação que contrasta com a realidade de muitos países europeus, onde persiste resistência a utilização de telemóveis na sala de aula. Convêm observar que a utilização de smartphones com as aplicações do GeoGebra permite ao estudante ter um elemento tecnológico de excelência para a aprendizagem da matemática e das ciências sem custos acrescidos, evitando a necessidade de utilização de outros dispositivos como sejam máquinas de calcular científicas ou gráficas, computadores pessoais no trabalho quotidiano inerente ao estudo destas disciplinas. Estes dispositivos podem utilizar o GeoGebra Exam, aplicação que permite com segurança ou seu uso em exames, já usada em alguns países da Europa e nos Estados Unidos.

Considerando o carácter multifacetado do GeoGebra que leva à que este software seja considerado por vários investigadores, a nível global e nacional, como um ambiente de aprendizagem matemática privilegiado (Dos Santos & Trocado, 2016; Santos & Quaresma, 2013; Quaresma, Santos & Bouallegue, 2013; Stahl, Rosé, O'Hara, & Powell, 2010). De facto, o franco desenvolvimento e expansão global das plataformas GeoGebra Resources, disponibilizando recursos livres para estudantes e professores, e GeoGebra Groups, com um sistema de gestão de aprendizagem que permite a estão da aprendizagem de grupos de estudantes com fedback personalizado. Estas plataformas continuam a ser estudadas e melhoradas pelo GeoGebra Team com a colaboração da comunidade global de educadores matemáticos que usam quotidianamente o software, nomeadamente com a evolução de sistemas de feedback automático que estão em desenvolvimento graças as ferramentas de raciocínio automatizado existentes no GeoGebra (Abánades, Botana, Kovács, Recio, & Sólyom-Gecse, 2016; Kovács, Recio, & Velez, 2019; Dana-Picard, 2019; Quaresma, Santos, Graziani, & Baeta, 2019).

Para alem da utilização do GeoGebra no ensino e aprendizagem da matemática, este software está a ser muito utilizado como recurso no ensino das ciências em geral, há múltiplos exemplos de utilização na área da física (Budak, Sahin & Dogan, 2018; Marciuc & Miron, 2018; Costa, 2018; Scheuber, 2018), da modelação de problemas de biologia (Jagiełło, Walczak, Iszkuło, Karolewski, Baraniak, & Giertych, 2019), inclusive nas artes que pelas suas capacidades CAD quer pela sua plasticidade e aplicação nas artes eletrónicas. Numa visão atual e para o futuro do ensino e aprendizagem da Matemática, onde o conhecimento e o pensamento matemático é um mediador das experiencias positivas em ciência, engenharia, tecnologia e as artes - STEAM, o GeoGebra tem sido usado em vários projetos onde participam estudantes. Professores e investigadores de vários países (Fenyvesi, Park, Choi, Song & Ahn, 2016; Budinski, Lavicza & Fenyvesi, 2018; Lavicza, Fenyvesi, Lieban, Hohenwarter, Mantecon & Prodromou, 2018).

A versão do GeoGebra usada e testada na realização das tarefas que apresentamos foi a versão GeoGebra Classic para Desktop, disponíveis nos seguintes links.

Windows: https://download.geogebra.org/package/win
MAC: https://download.geogebra.org/package/mac
Para outras versões, ver o fim da página do link : https://wiki.geogebra.org/en/Reference:GeoGebra_Installation

## *Enquadramento das tarefas*

As tarefas aqui apresentadas fazem parte um estudo mais amplo, que se desenvolve, considerando os ambientes tecnológicos: Geogebra e STEAM, como suportes da formação inicial e contínua de professores e como opurtunidade de desenvolver processos de reflexão e de alteração sustentada das praticas associadas ao ensino e aprendizagem da matemática.

Cada uma das tarefas esta pensada no sentido de poder ser associada a realização de um projeto a desenvolver no espaço escolar, na disciplina de matemática bem como outras disciplinas que se podem associar ao projeto. Na primeira versão, os guiões das tarefas que aqui se propõe desenvolvem conteúdos matemáticos, contudo posteriormente estes protocolos serão estendidos de modo que a sua utilização se alargue a professores das áreas das ciências e das artes. Pretende-se que estas propostas sejam pólos aglutinadores de conhecimentos, contemplando no trabalho escolar varias estilos de aprendizagem, desenvolvendo aprendizagens e competências nos alunos.

Nesta primeira fase foram elaboradas nove tarefas a saber:

*Tarefa 1 – Modelar poligonais fechadas sobre um cubo a a partir de uma planificação*

*Tarefa 2 – Famílias de funções reais de variável real*

*Tarefa 3 – Modelação de uma Reta, de um Segmento de Reta e de uma Semirreta, como representação cartesiana do gráfico de uma função real de variável real*

*Tarefa 4 – Modelação de superfícies contidas em planos e definidas por funções bivariadas.*

*Tarefa 5 – Modelação de sólidos de revolução, uso de splines*

*Tarefa 6 – Modelação de conchas*

*Tarefa 7 – Modelação por Curvas de Nível*

*Tarefa 8 – Modelação com Epiciclos*

*Tarefa 9 – Modelação de Airfoils com Números Complexos*

Estas tarefas foram concebidas por José Manuel Dos Santos, Astrigilda Silveira e Alexandre Trocado; sendo os dois primeiros formadores, para os PALOP. As tarefas só podem ser usadas em momentos de formação de professores e/ou adaptadas para serem usadas em sala de aula. O uso das tarefas está sujeito a reporte para os seguintes endereços electrónicos: santosdossantos@ese.ipp.pt e astrigilda.silveira@docente.unicv.edu.cv. Desde já, agradecemos a colaboração de todos os interessados em usar estas tarefas, através do envio de *feedback* ou comentários sobre o seu uso, esperando que as mesmas sejam úteis a professores e alunos no ensino e aprendizagem da Matemática.

## *Tarefa 1 – Modelar poligonais fechadas sobre um cubo a partir de uma planificação*

Pretende-se criar um cubo, podendo variar a medida da aresta com um selector. Iniciemos o **GeoGebra**, dirija-se ao menu Vista, de modo a configurar a sua janela de trabalho visualizando a **Folha Algébrica**, a **Folha Gráfica 2D** e a **Folha Gráfica 3D** .

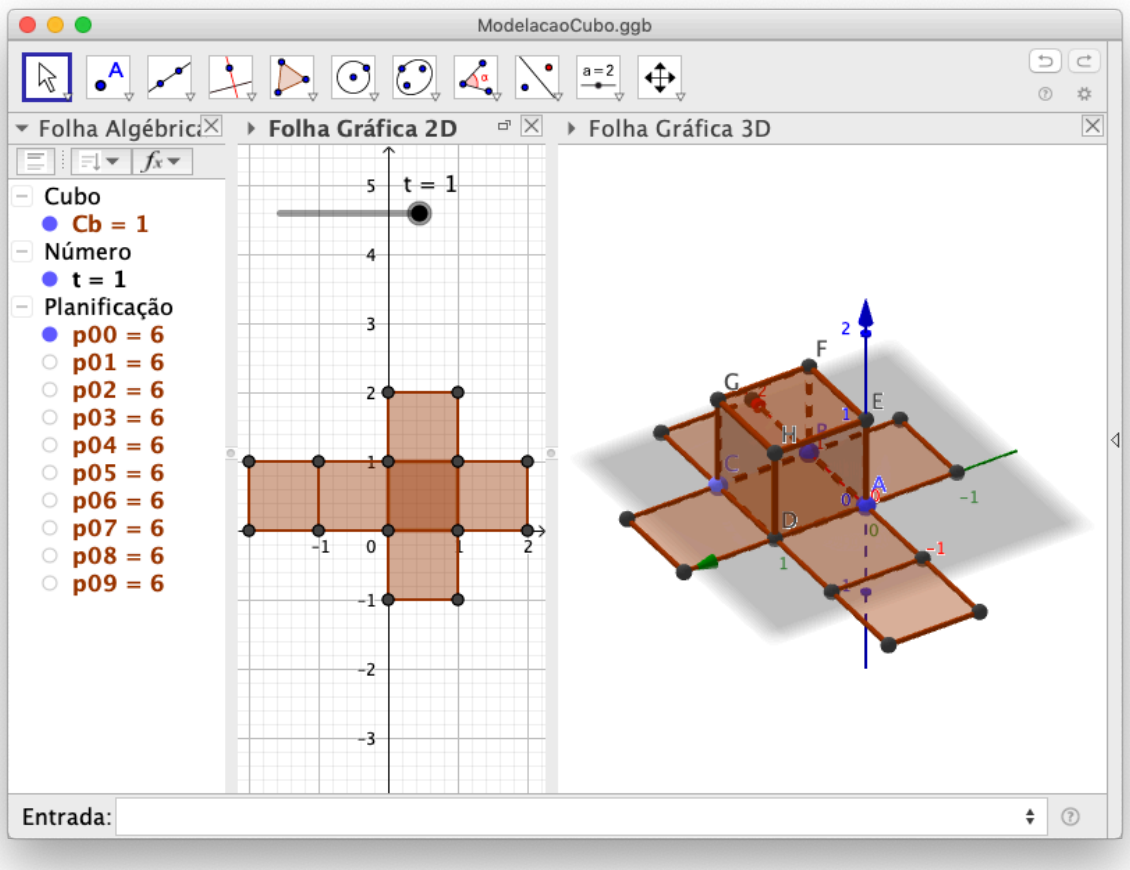

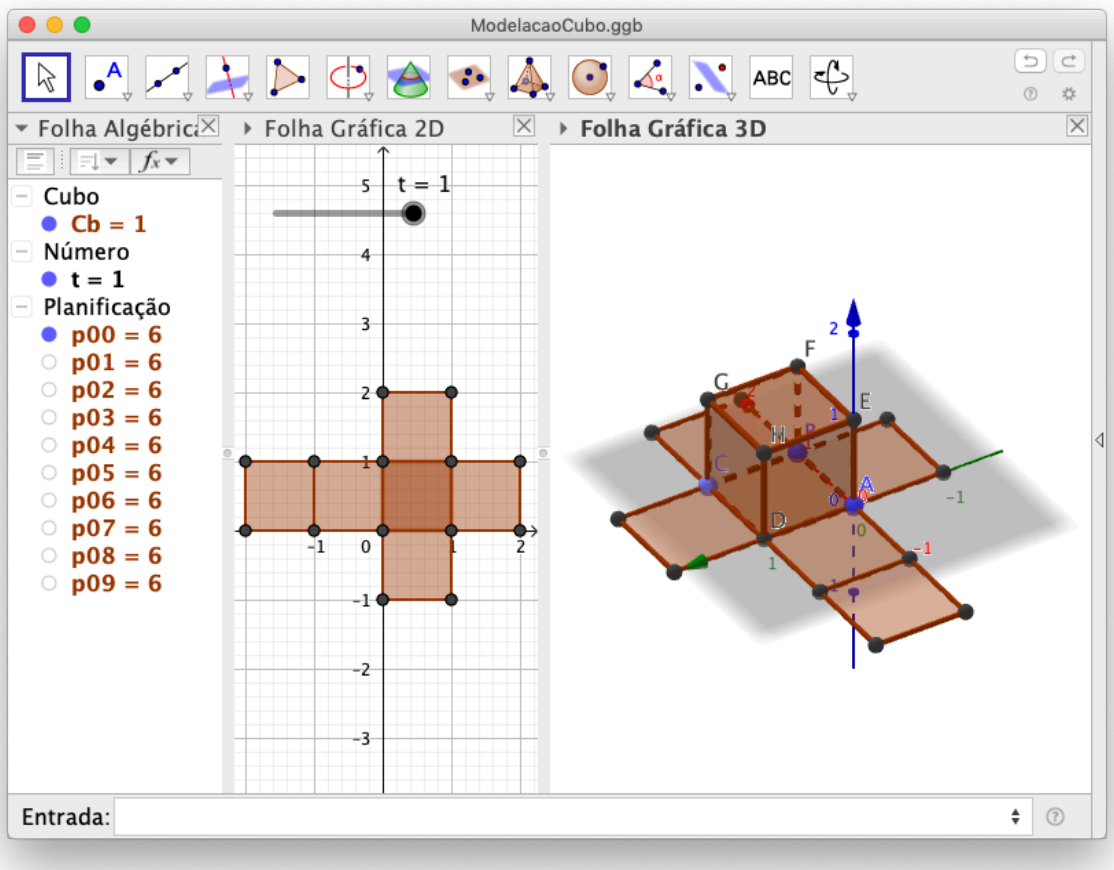

1. Com a **Folha Gráfica 3D** ativa, introduza no campo de **Entrada**, os pontos A e B, dando entrada dos comandos:
   a) `A=(0,0,0)`
   b) `B=(1,0,0)`
2. Com a **Folha Gráfica 3D** ativa, introduza no campo de **Entrada**, o comando que lhe vai permitir criar o cubo:
   a) `Cb=Cubo(A,B)`
   b) Pode em alternativa usar a ferramenta 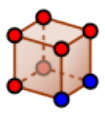, cubo, e selecionando os pontos A e B com o cursor.

   Note para selcionar pontos tem de ter ativa a ferrameta 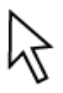 , mover.
3. Com a **Folha Gráfica 2D** ativa, introduza no campo de **Entrada**, o comando que lhe vai permitir criar o seletor que vai revelar a planificação do cubo:
   a) `t=Seletor(0, 1, .1, 1, 100, false,true,false, false)`
4. A planificação "standard" do cubo pode ser obtida, com a **Folha Gráfica 3D** ativa:
   a) dando entrada no campo de Entrada de `p00=Planificação[Cb, t]`.
   b) ou usando a ferramenta 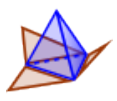 planificação.
5. Desloque o seletor `t` e verifique o desenvolvimento do desdobrar da superfície do cubo.
6. A partir de comandos, o GeoGebra permite obter outras planificações do cubo. Por exemplo, use algum dos seguintes comandos:
   a) `p01=Planificação(Cb, t, faceABCD, arestaBC, arestaCD, arestaCG, arestaDH)`
   b) `p02=Planificação(Cb, t, faceABCD,arestaAD,arestaAE,arestaBF,arestaCD)`
   c) `p03=Planificação(Cb, t, faceABCD, arestaAD,arestaBC,arestaBF,arestaCG)`
   d) `p04=Planificação(Cb, t, faceABCD, arestaAD,arestaCD,arestaEH,arestaGH)`
   e) `p05=Planificação(Cb, t, faceABCD, arestaAD,arestaBC,arestaCG,arestaDH)`
   f) `p06=Planificação(Cb, t, faceABCD, arestaAE,arestaBC,arestaBF,arestaCD)`
   g) `p07=Planificação(Cb, t, faceABCD, arestaAB,arestaAE,arestaCD,arestaDH)`
   h) `p08=Planificação(Cb, t, faceABCD, arestaAB,arestaBF,arestaCD,arestaDH)`
   i) `p09=Planificação(Cb, t, faceABCD, arestaAB,arestaAE,arestaBC,arestaBF)`
7. Numa nova aplicação do GeoGebra repita os pontos de 1 a 5 e:
   a) faça com que o seletor t tome o valor um;
   b) usando a ferramenta 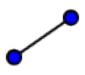 , segmento de reta, ou em alternativa a ferramenta 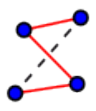 , linha poligonal, marque uma linha poligonal com três segmentos de reta.
   c) faça com que o seletor `t` tome o valor zero e observe;
   d) ajuste ou marque outra linha poligonal de modo esta seja uma linha poligonal fechada quando `t=0`.

   Note que o comando planificação só se aplica a poliedros platónicos.

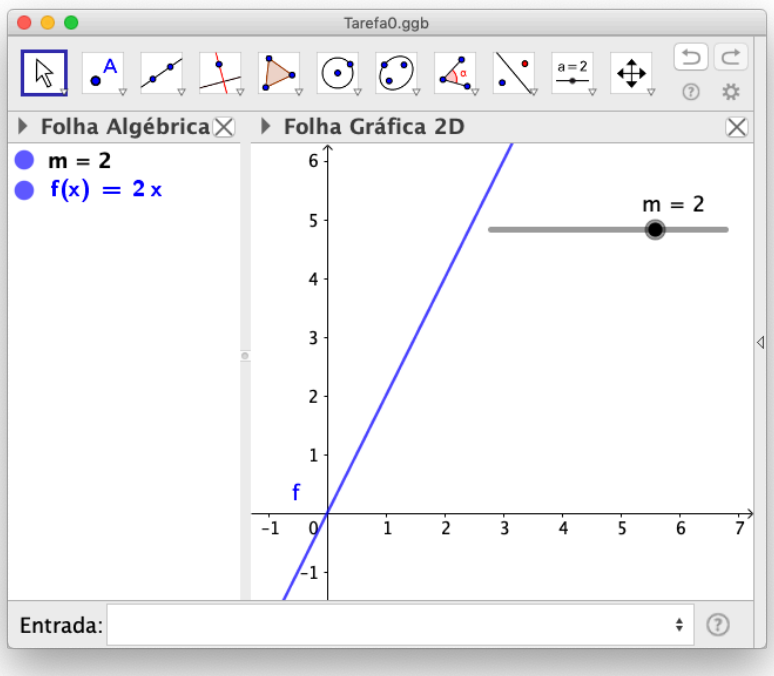

Imagem de aplicação do GeoGebra: na ***Folha Gráfica 2D*** está a representação cartesiana do gráfico de uma função; o parâmetro *m* e a função f está definida como *f(x)=m x* na *Folha Algébrica*.

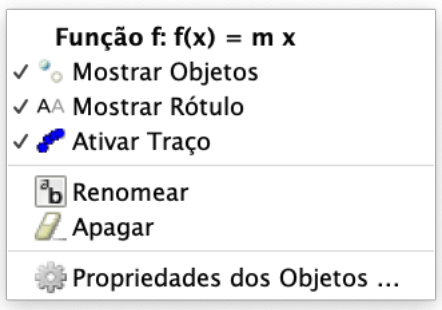

*A saber:*

- Ferramenta seletor, a=2, e usar a campo de **Entrada** para o obter.
  Um ponto parâmetro pode ser representado por um seletor no GeoGebra. Os seletores também podem ser obtidos usando o campo de **Entrada** com o comando
  `Seletor( <Min>, <Max>, <Incremento>, <Velocidade>, <Largura>, <Ângulo>, <Horizontal>, <Animação>, <Aleatório> )`

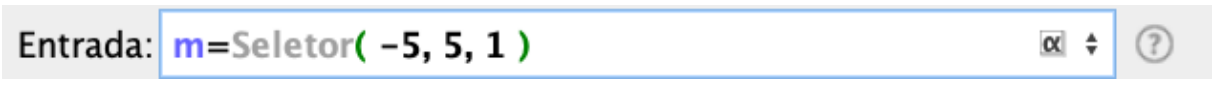

- Funções
  Para obter uma função usamos a escrita habitual no campo de **Entrada**.

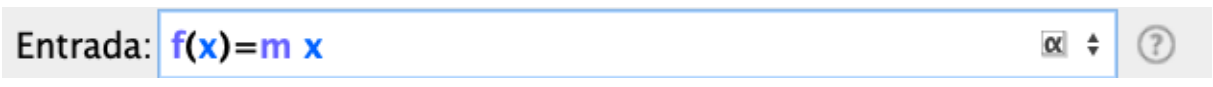

  Por exemplo, sendo m um número ou um parâmetro, escrevendo `f(x)=m x` no campo de Entrada, definimos a função identidade de $\mathbb{R}$ em $\mathbb{R}$, obtendo a representação cartesiana do gráfico de *f* numa das folha gráficas bidimensionais.
- Comando sequência.
  Para representar vários elementos da família de funções podemos usar no campo de **Entrada** o comando sequência. Por exemplo, `Ff=Sequência(m x,m,-5,5,1)`. Em alternativa podemos ativar o traço de *f*, , acessível pressionando o botão direito do rato sobre *f*.

1. Obtenha no GeoGebra representações de:

   1.1. $f(x)=ax^2, a\in\{-2,\dots,2\}$;

   1.1. $g(x)=mx+b, m\in\{-2,0,2\}\wedge b\in\{-2,0,2\}$;

   1.3. $h(x)=sqrt(x-a)+b$ onde $a,b\in\{-1,1\}$.

2. Considere a família $l(x)=ax^2+bx+c$, com $a,b,c\in[-5,5]$.

   2.1. Represente os parâmetros *a*, *b* e *c* por seletores com incremento de um.

   2.2. Pressione o botão direito do rato sobre *l*, aceda a ***Propriedades dos objetos*** e no separador ***Avançado***, proceda como na imagem.

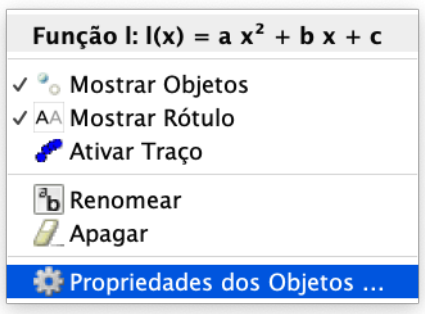

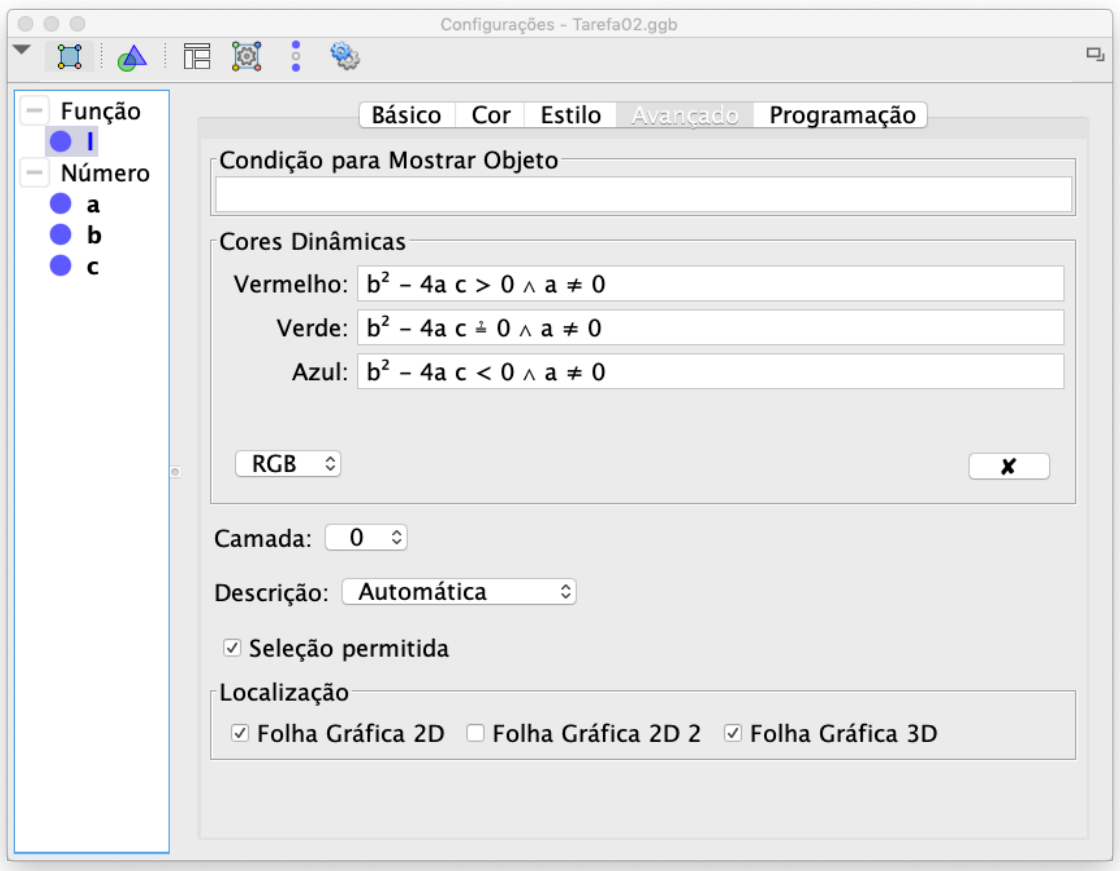

Antes de voltar a aplicação do GeoGebra, e movimentar os seletores o que espera que aconteça em relação as cores da representação cartesiana do gráfico?

Confirme as suas conjeturas!

3. Adapte a aplicação construída em 2 para estudar a representação cartesiana do gráfico da família de funções reais de variável real definidas por $m(x)=a|x-b|+c; a,b,c\in\mathbb{R}$.

## *Tarefa 3 – Modelação de uma Reta, de um Segmento de reta e uma Semirreta, como representação cartesiana do gráfico de uma função real de variável real*

***A saber sobre pontos e funções no GeoGebra.***

Ferramenta ponto, , e usar o campo de entrada para o marcar ou obter as suas coordenadas.

Um ponto *P* pode ser obtido no GeoGebra com a ferramenta Ponto, sobre as folhas gráficas, as suas coordenadas correspondem ao comando *(x(P),y(P))*, em 2D, ou a *(x(P),y(P),z(P))* em 3D.

Os pontos, também podem ser obtidos usando a barra ou campo de **Entrada**.

Entrada: B=(3,4)

Por exemplo, escrevendo B=(3,4), se apenas escrevermos (3,4) o software atribui ao ponto um nome que corresponderá a primeira letra maiúscula, ou sequencia de letras maiúsculas disponíveis.

Funções

Para obter uma função usamos a escrita habitual no campo de entrada. Por exemplo, escrevendo *f(x)=x* no campo de Entrada, definimos a função identidade de $\mathbb{R}$ em $\mathbb{R}$, e obtemos a representação cartesiana do gráfico de *f* numa das folha gráficas bidimensionais.

O comando Função( <Função>,<Valor de x Inicial>,<Valor de x Final>, permite obter uma função real de variável real cujo domínio é uma parte de $\mathbb{R}$. Por exemplo o comando g(x)=Função(x^2,0,+∞), define uma função quadrática restrita ao intervalo .

1. Marque, na **Folha Gráfica 2D**, com ajuda da ferramenta ponto, , um ponto *A*.

2. Use o campo de **Entrada** para obter o ponto *B (3,4)*. Acautele que os pontos *A* e *B* não coincidem.

3. Represente, na **Folha Gráfica 2D**, a reta *AB* sem usar a ferramenta , reta, ou o comando reta.

4. Defina uma função cuja representação cartesiana do seu gráfico seja [*BC*] onde C(4,3). Mais uma vez não use usar a ferramenta , segmento de reta, nem o comando associado

5. Considere uma função cuja representação cartesiana do seu gráfico é a semireta *CD*, sem usar a ferramenta , semirreta. Em que condições tal função existe?

6. Sejam *E* e *F* dois pontos cuja ordenada é igual a três. Sem usar a ferramenta , ponto médio, nem a ferramenta , semicírculo, defina uma função cuja representação cartesiana do seu gráfico seja um semicírculo.

Note que os comandos para obter no campo de Entrada a raiz cúbica e a raiz quadrada são cbrt(<x>) e sqrt(<x>), respectivamente,

6.1. Nestas condições quantas funções podem ser definidas?

6.2. Se os pontos *E* e *F* não tivessem a mesma ordenada seria possível ainda definir uma função cuja representação cartesiana do seu gráfico estivesse contida num semicírculo?

Funções de $\mathbb{R}^2$ em $\mathbb{R}$.

O GeoGebra também reconhece funções de $\mathbb{R}^2$ em $\mathbb{R}$. A representação cartesiana do gráfico de uma função de $\mathbb{R}^2$ em $\mathbb{R}$ define um subconjunto de pontos em $\mathbb{R}^3$ ao qual esta associado uma superfície. Por exemplo, `Função(<Expressão>, <Variável 1>, <Valor Inicial>, <Valor Final>, <Variável 2>,<Valor Inicial>,<Valor Final>))`, permite obter na **Folha Gráfica 3D** a representação cartesiana do gráfico de uma função bivariada, que corresponde a uma superfície em $\mathbb{R}^3$ .

No caso da imagem seguinte, para a função *h(x,y)* pedida em 1.a.i., corresponde a um losango contido no plano de equação cartesiana *z=x+y*.

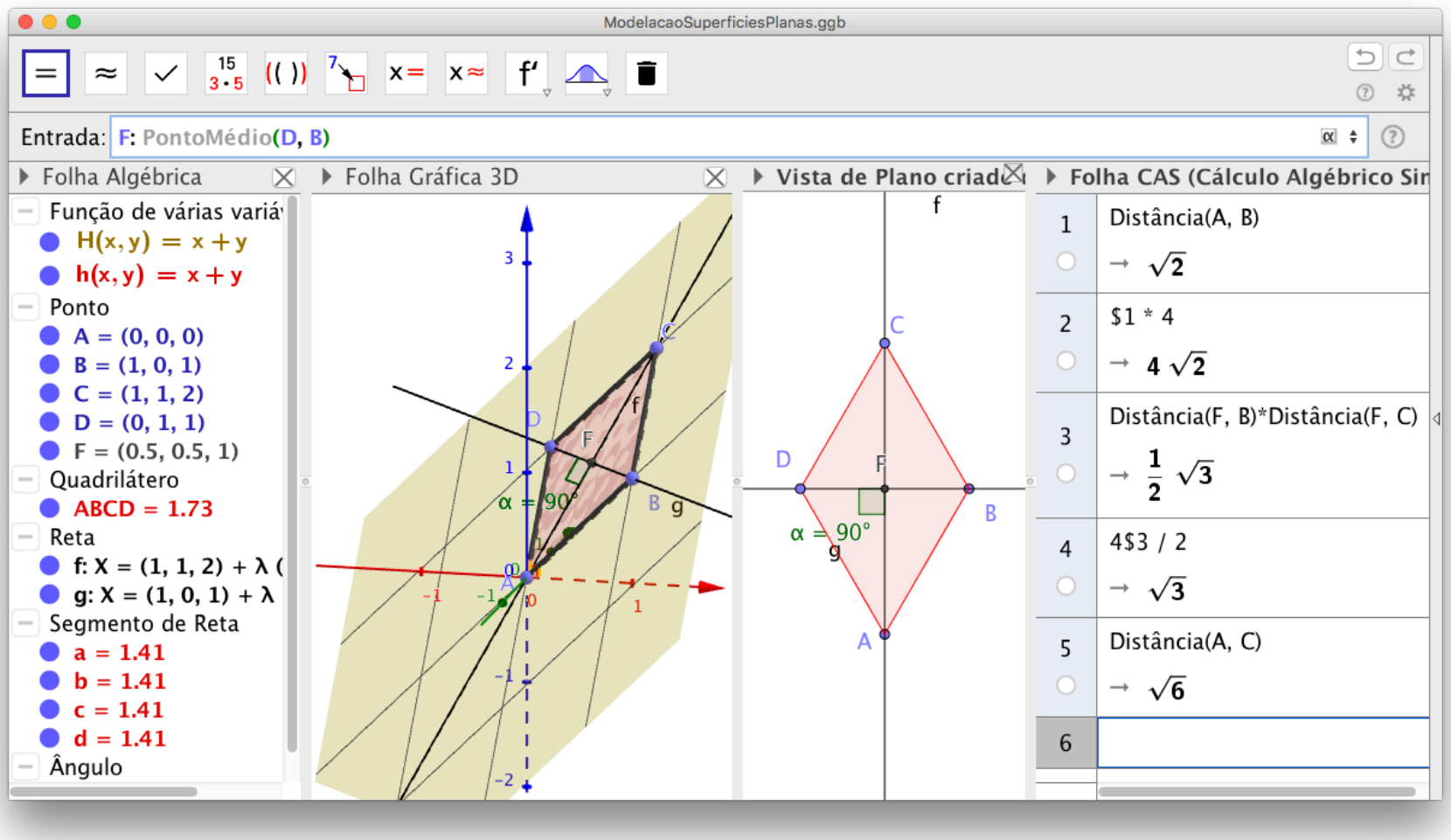

1. Seja a função H definida por *H(x,y)=x+y*, onde *x* e *y* são números reais.
   a. No campo de **Entrada** e use o comando:
      i. `h(x,y)=Função( x+y,x,0,1,y,0,1)` para representar a restrição da função *H* restrita ao conjunto $[0,1]x[0,1]\subset\mathbb{R}^2$;
      ii. `H: z=h(x,y)`
   b. Observe a **Folha Gráfica 3D** e compare as representações obtidas
2. Represente:
   a. os pontos de coordenadas:
      i. `A: (0,0,h(0,0))` ;
      ii. `B: (1,0,h(1,0))` ;
      iii. `C: (1,1,h(1,1))` ;
      iv. `D: (0,1,h(0,1))` ;
   b. o polígono [*ABCD*] com o comando `Polígono(A,B,C,D)`.

   Note que poderia também utilizar a ferramenta 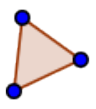, polígono.

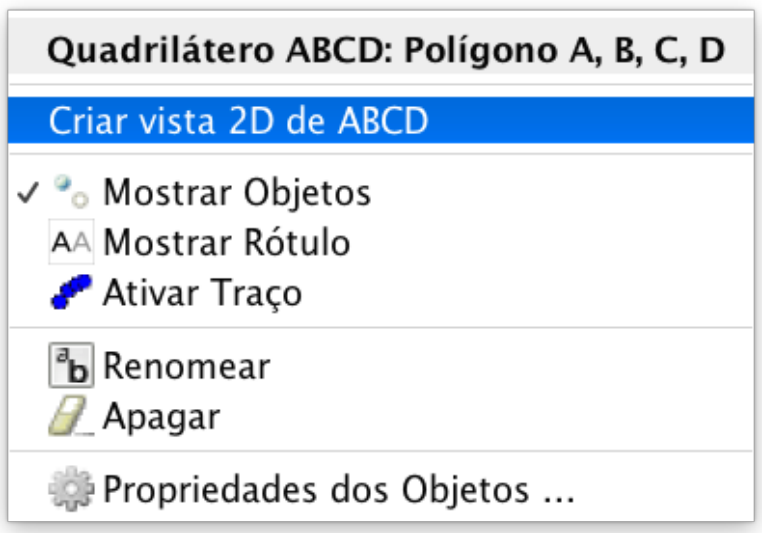

3. Obtenha a **vista 2D de [ABCD]** e com as ferramentas da geometria bidimensional verifique que o polígono é um losango.
4. Use a **Folha CAS** para obter valores exatos de algumas medidas associadas a [*ABCD*].
5. Considere a função *G* definida por *G(x,y)=|x+y|*, onde *x* e *y* são números reais.
   a. No campo de **Entrada**, use o comando *G: z=abs(x+y)* para obter a representação do gráfico cartesiano de *G* **na Folha Gráfica 3D**.
   b. Obtenha a representação cartesiana de uma restrição de *G* que não corresponda a uma superfície plana.

## *Tarefa 5 – Modelação de sólidos de revolução, uso de splines*

Para um objeto, cujo modelo seja um sólido de revolução, poderia tirar-se uma foto em perfil e proceder como abaixo se refere. Por outro lado, recorrendo à aplicação GeoGebra AR (realidade aumentada), este processo poderia decorrer em vídeo e em tempo real. Aqui deixamos uma versão simplificada desta tarefa.

1. Encontre uma imagem de um sólido de revolução de perfil na Web. Insira a imagem na **Folha Gráfica 2D**, com a ferramenta inserir imagem, .

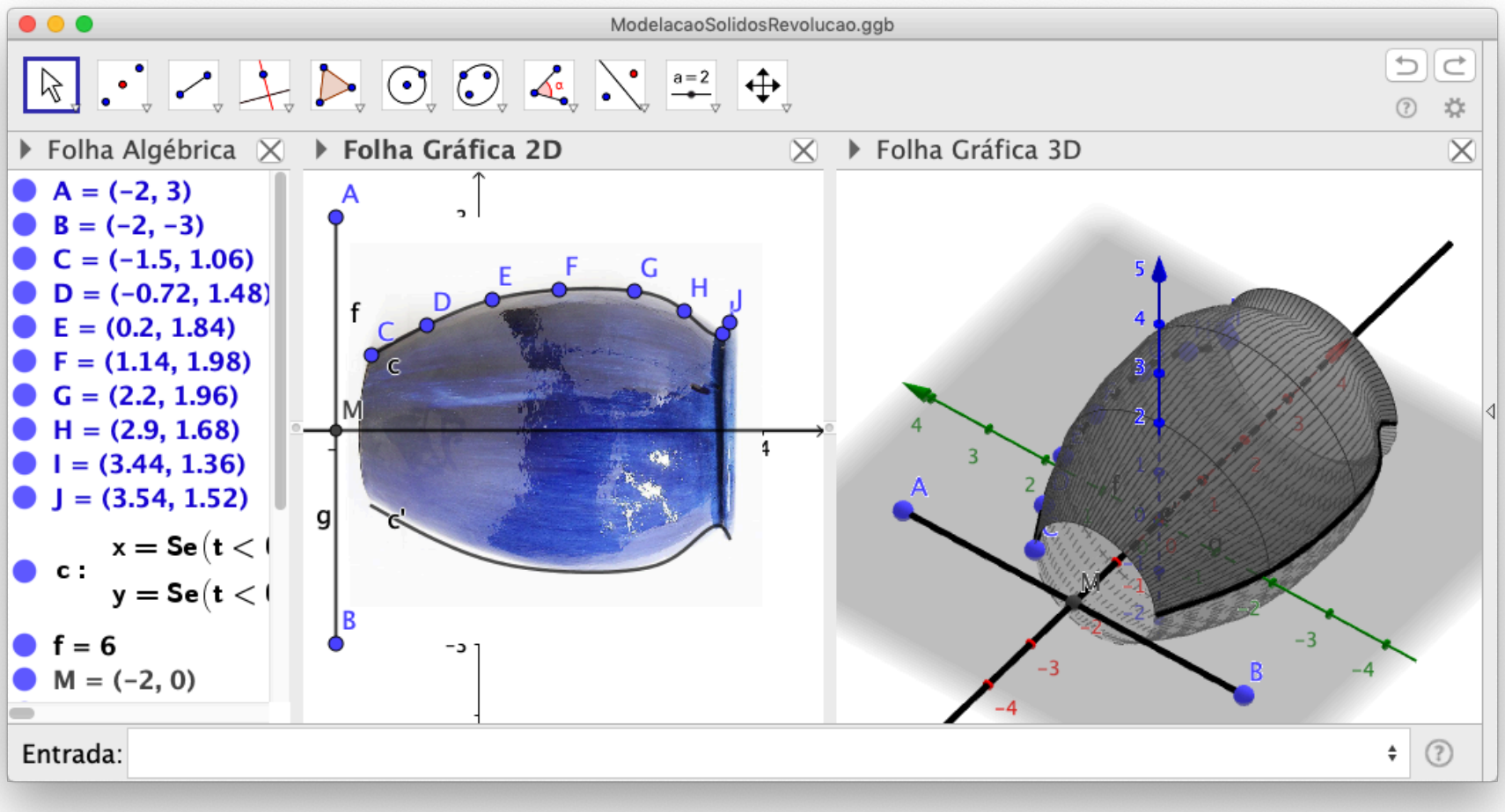

2. Posicione a imagem, com a ferramenta mover, , de modo que o seu eixo de autoreflexão coincida com um dos eixos *Ox ou Oy*.

    a. Marque o ponto médio de [*AB*] , `M=PontoMédio(A,B)`, ou use .

    b. Represente o segmento de reta [*AB*], usando a ferramenta , segmento de reta, seleccionando os pontos *A* e *B*.

    c. Pelo ponto médio de [*AB*] marque a perpendicular a [*AB*] em *M* com a ferramenta .

    Observe que: podia ter determinado a reta anterior com a ferramenta , mediatriz; os eixos do referencial são designados por `EixoOx`, `EixoOy` e `EixoOz`, sempre que necessitar de os usar nos comandos.

3. Com a ferramenta ponto, , distribua pontos em parte do contorno da imagem que seja invariante por reflexão em relação a uma reta que possa ser identificada ou traçada na imagem.
4. Use o comando `Spline(<Lista de Pontos>,<Ordem ≥ 3>)`, substituindo a lista de pontos pelas letras maiúsculas que representam os pontos marcados no ponto anterior, separados por virgulas, e colocados entre chavetas.

    No caso da imagem acima o comando `c=Spline({C,D,E,F,G,H,I,J},3).`

5. Obtenha a reflexão da curva `c` na perpendicular obtida em 2.c. com a ferramenta e ajuste a curva `c`.
6. Ative a **Folha Gráfica 3D**, use o comando superfície lateral no campo de Entrada para modelar tridimensionalmente o objeto.

    No imagem acima foi usado o comando `Sc=SuperfícieLateral(c,2π,EixoOx)` no campo de **Entrada**, para modelar o objeto tridimensional.

7. Manipule a vista tridimensional do objeto usando a ferramenta, , rodar a vista 3D.

## *Tarefa 6 - Modelação de conchas*

O GeoGebra permite representações tridimensionais que podem ser usados para modelar objetos. Esta capacidade pode ser usada, por exemplo, para modelar conchas de moluscos. Há vários estudos sobre modelação de conchas que podem ser dados como bibliografia aos alunos (Cortie, 1989; Pappas & Miller, 2013) para a realização de projetos, estabelecendo relações com outras ciências, e usando o software GeoGebra.

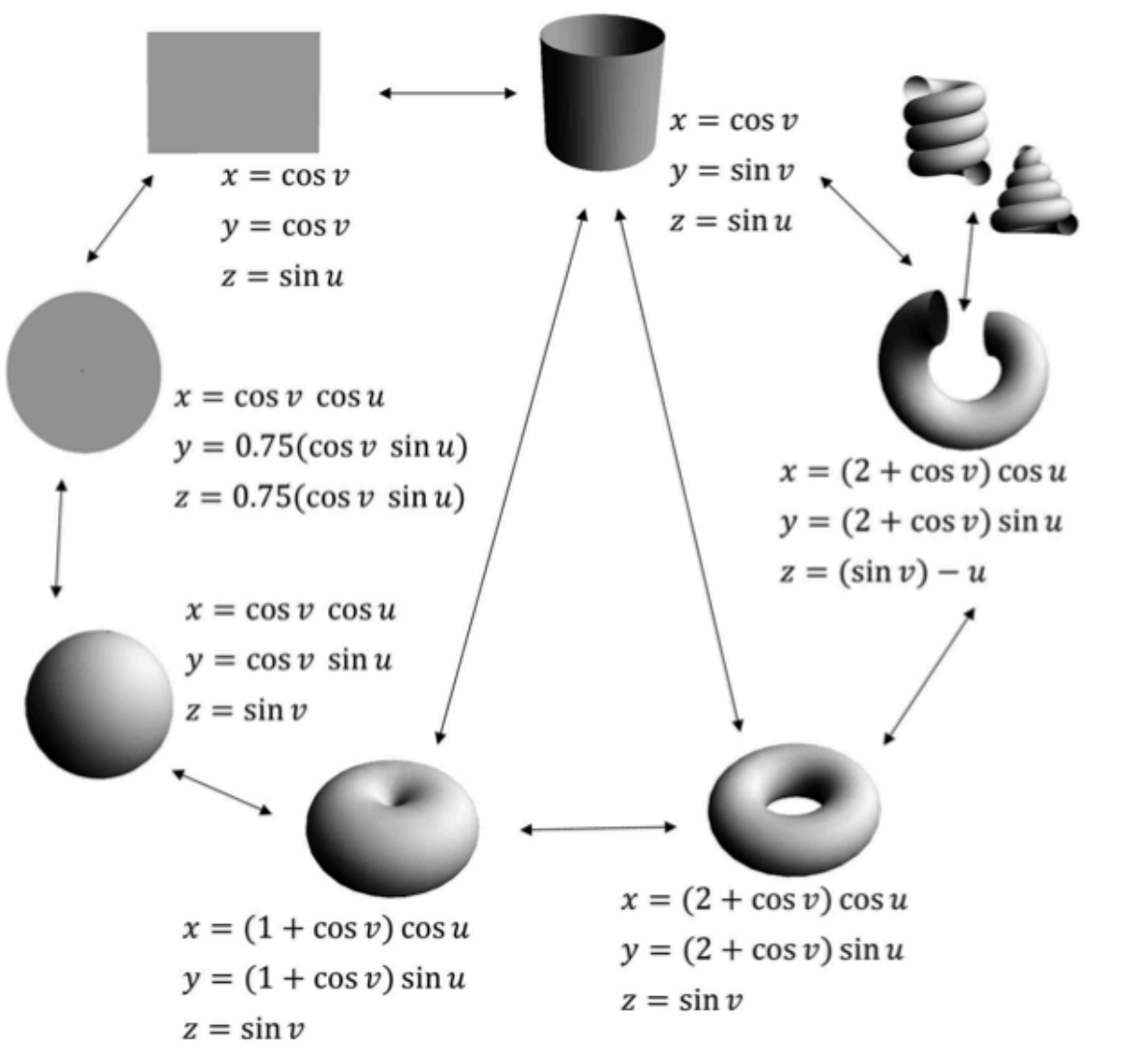

Figura 1a. Das equações paramétricas de superfícies geométricas básicas até a modelação de conchas (Pappas & Miller, 2013). doi:10.1371/journal.pone.0077551.g001

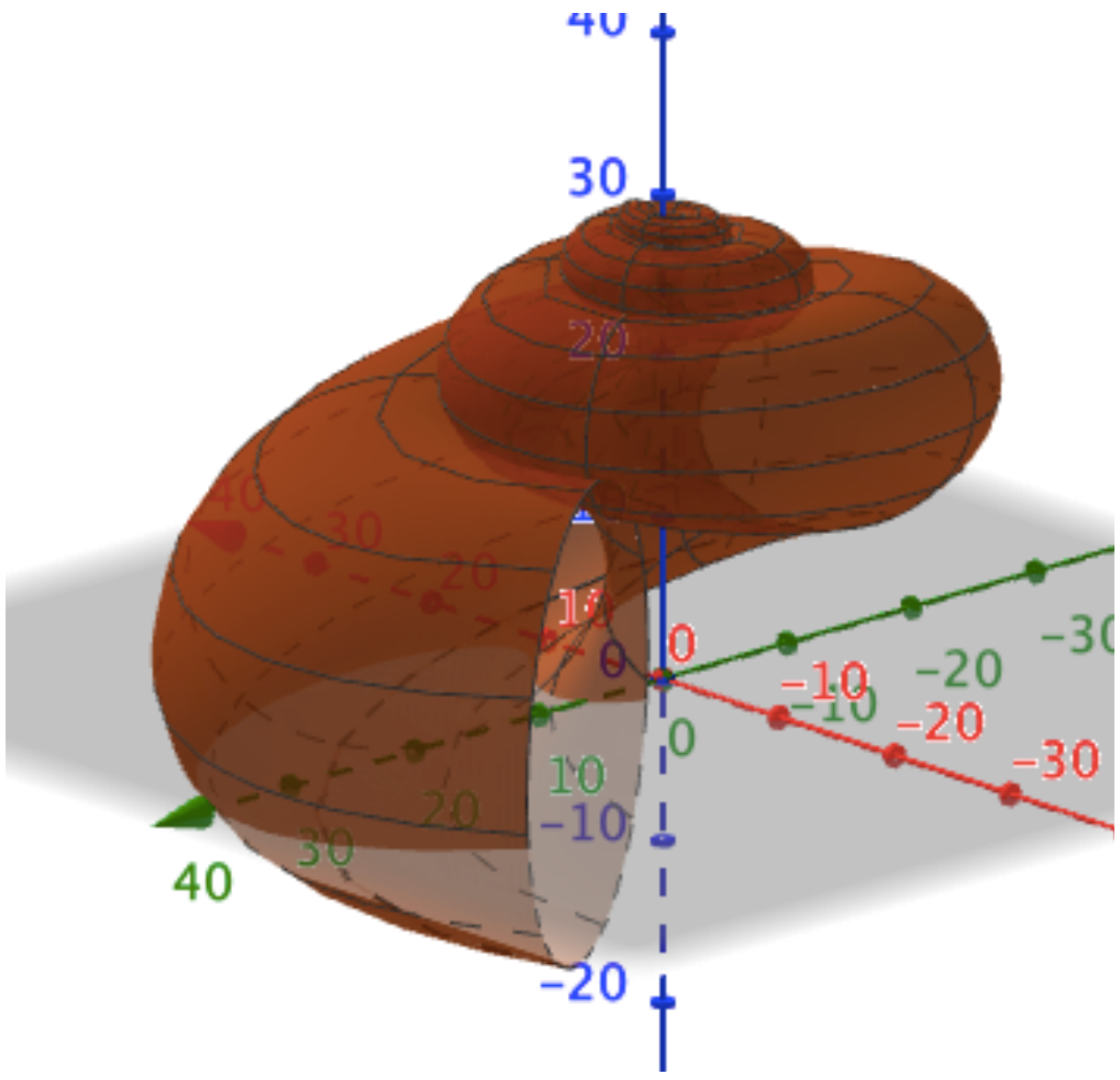

Figura 2a. Modelação no GeoGebra obtida com o comando Superfície Lateral

1. Modele a superfície de uma esfera e um toro, usando o comando Superfície lateral:

```
SuperfícieLateral( <Expressão>, <Expressão>, <Expressão>, <Variável 1>, <Valor Inicial>, <Valor Final>, <Variavel 2>, <Valor Inicial>, <Valor Final> )
```

Sugestão: atenda a informação da Figura 1a.

2. Para modelar a imagem da Figura 2a. no GeoGebra, comece por introduzir os comandos para definir: i) os parâmetros da “Concha”; ii) as funções componentes da parametrização; iii) os limites do domínio da parametrização; iv) a superfície.

**i)**

```
r=Seletor(0,1,1/100,1,100)
f=Seletor(0,1,1/100,1,100)
a=Seletor(0,5,1,1,100)
b=Seletor(0,5,1,1,100)
c=Seletor(0,5,1,1,100)
d=Seletor(0,5,1,1,100)
e=Seletor(0,5,1,1,100)
h=Seletor(-
30,30,1/100,1,100)
```

**ii)**

```
X(u,v)=r*exp(1)^(f*u)*(d+a*cos(v))*cos(c*u)
Y(u,v)=r e^(f u) (d +a cos(v))sin(c u)
Z(u,v)=r e^(f u) (-1.4e+b sin(v))+h
```

**iii)**

```
u1=Seletor(-5,5,1/100,1,100)
u2=Seletor(-5,5,1/100,1,100)
v1=Seletor(-5,5,1/100,1,100)
v2=Seletor(-5,5,1/100,1,100)
```

**iv)**

```
s=SuperfícieLateral(X(u,v),Y(u,v),Z(u,v),u,u1,u2,v,v1,v2)
```

## *Tarefa 7 - Modelação por Curvas de Nível*

Uma curva de nível corresponde a representação de uma linha imaginária que une todos os pontos de igual "altitude", com a mesma cota, de uma mesma superfície. As curvas de nível correspondem a representação da intersecção de superfícies com planos paralelos ao plano *xOy*.
As superfícies podem ser obtidas no GeoGebra a partir de uma equação nas variáveis reais *x*, *y* e *z* ou, em casos particulares, como uma função de $\mathbb{R}^2$, ou parte de $\mathbb{R}^2$, em $\mathbb{R}^3$. Note que o GeoGebra permite representar na **Folha Gráfica 3D** a imagem de uma função de variável, ou variáveis reais, cujo conjunto de chegada é $\mathbb{R}^3$.

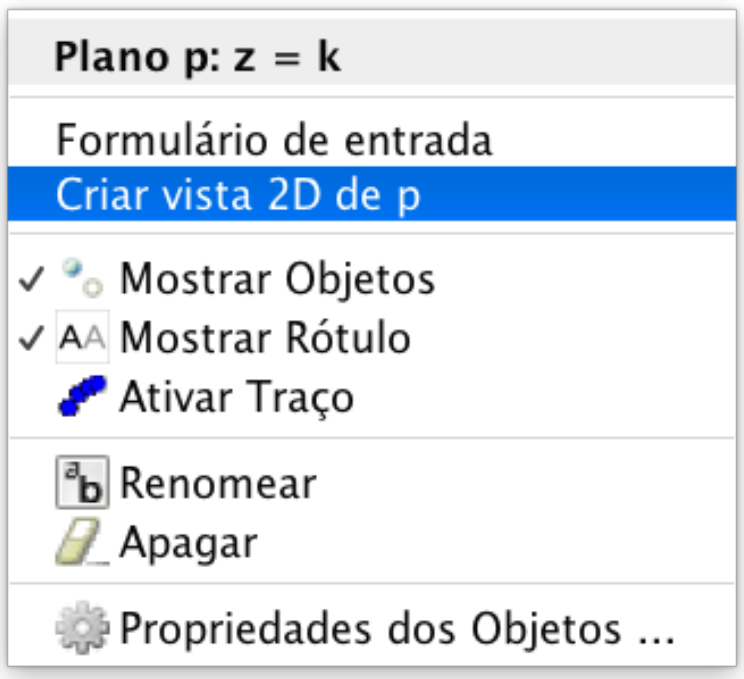

A representação da intersecção de uma superfície com os planos do tipo $z=k,\ k\in\mathbb{R}$ , pode ser obtida na **Folha Gráfica 2D**, **Folha Gráfica 2D 2**, ou na representação bidimensional fornecida na vistas auxiliares.
A intersecção de um plano com uma superfície pode ser obtida com a ferramenta 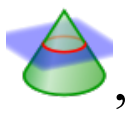, intersecção de duas superfícies ou com vários comandos associados, por exemplo,
`Interseção(<Objeto>,<objeto>)` .
A vista da curva no plano da curva pode ser obtida, clicando com o botão direito do rato no plano e seleccionando **Criar vista 2D de ...**

Nota: Um plano é uma superfície, poderia entender-se que o GeoGebra representaria a intersecção de duas quaisquer superfícies, ora este assunto não é um assunto fechado e ainda é alvo de estudo no GeoGebra (Dos Santos, 2017; Trocado, Gonzalez-Vega & Dos Santos, 2019), assim como se estudam algoritmos aplicáveis a outros softwares (Trocado, Gonzalez-Vega & Dos Santos, 2019).

1. Na imagem abaixo estão representadas curvas de nível. Nas três imagens obtidas no GeoGebra está a representação das curvas de nível obtidas pelos planos de equação $z=k$, $k\in\{0,1/4,1/2,2/3\}$.
   a. Quais são as superfícies que podem estar associadas a estas representações em curva de nível?
   b. Modele as situações envolvidas nas representações abaixo, em 2D e 3D, com o GeoGebra.

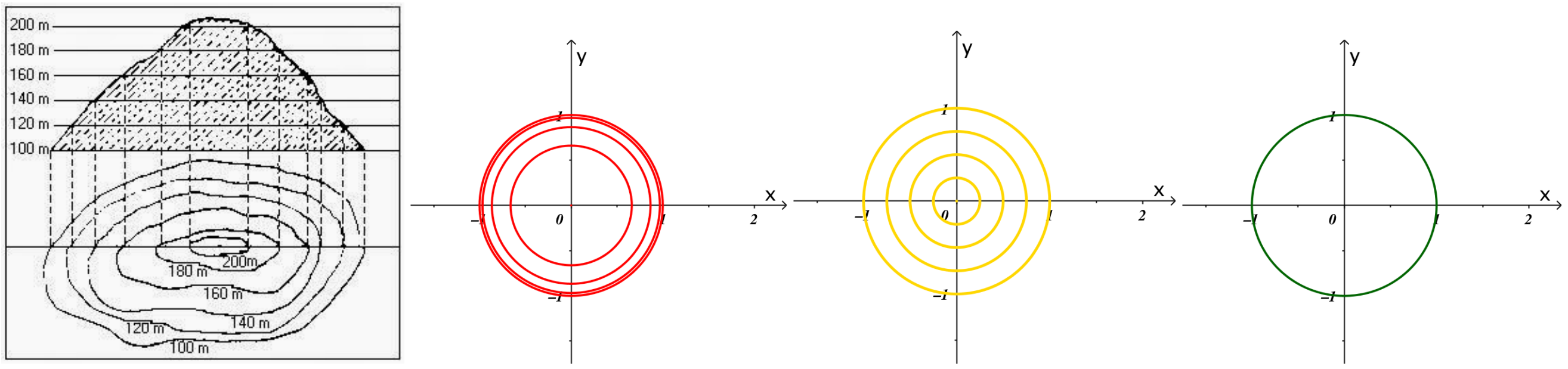

2. Na imagem abaixo estão representadas várias curvas de nível, determinadas por planos de equação $z=k$, na mesma janela de visualização, de uma superfície que se obtêm no GeoGebra a partir da representação cartesiana do gráfico de uma função, *M*, de $\mathbb{R}^2$ em $\mathbb{R}$.

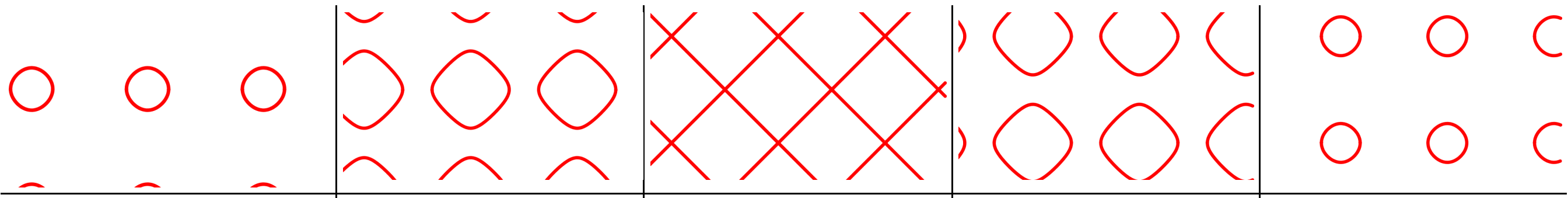

| k= - sqrt(2) | k= - sqrt(2)/2 | k=0 | k=sqrt(2) | k=sqrt(2) |
|---|---|---|---|---|

   a. Qual poderá ser a expressão analítica da função *M*?

Observação: Tarefa adaptada de atividades de um Curso de Formação de Formadores em GeoGebra para Cabo Verde (Fortes,2019)

3. Execute no GeoGebra a lista de comandos abaixo. Com a aplicação construída certifique as respostas dadas na questão anterior.

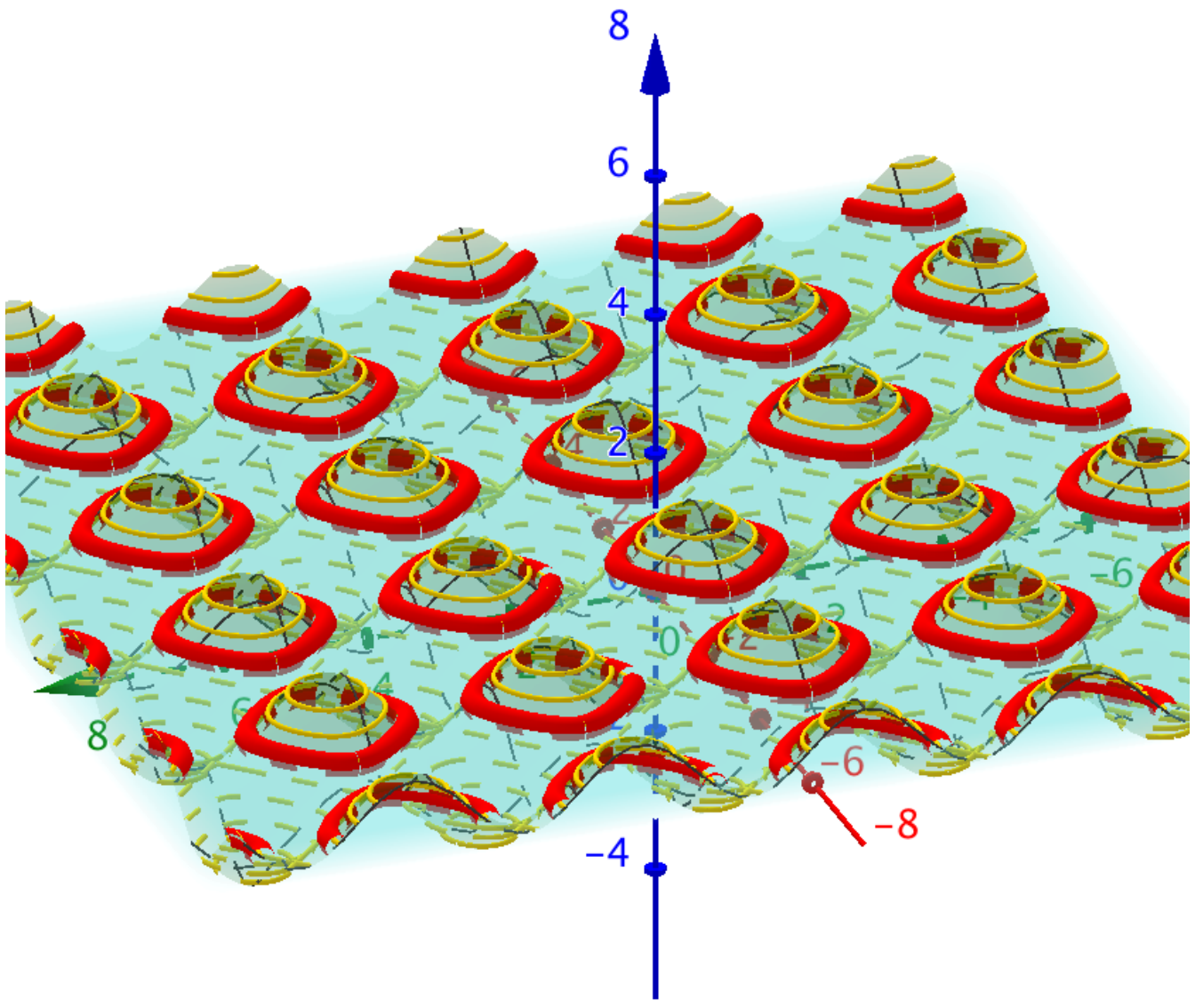

i) Para definir a superfície que corresponde ao Gráfico da função $f(x,y)=sin(x)^2+sin(y)^2$, $(x,y)\in\mathbb{R}^2$, na **Folha Gráfica 3D**, observar a representação ao dar entrada de:

```
M(x,y)=sin(x)^2+cos(y)^2
```

ii) Defina o plano de equação $z=k$, para cada k, com os comados:

```
k=Seletor(0, 2, 2/100, 1, 100)
p: z=k
```

(Faça com que o seletor apareça na **Folha Gráfica 2D**)

iii) Depois de **Criar a vista 2D de p**. Defina a curva de nível obtida da intersecção da representação cartesiana do gráfico de M e do plano p

```
iM=InterseçãoGeométrica(p, M)
```

iv) Para definir as curvas de nível usamos o comando:

```
SM=Sequência(InterseçãoGeométrica(z=k, M), k, 0, 2,1/50)
```

Nota: No caso do objeto M ser uma quádrica, obteríamos a projecção de todas as curvas de nível no plano xOy, usando os dois comandos seguintes:

`t=Seletor(0, 1, 1/100, 1, 100)`, colocando-o na ***Folha Gráfica 2D***.

```
PMt=Sequência((x(Ponto(Elemento(SM,j),t)),y(Ponto(Elemento(SM,j),t))),j,1, Comprimento(SM),1)
```

```
LPM=Sequência(Lugar_Geométrico(Elemento(PMt,j),t),j,1,Comprimento(SM),1)
```

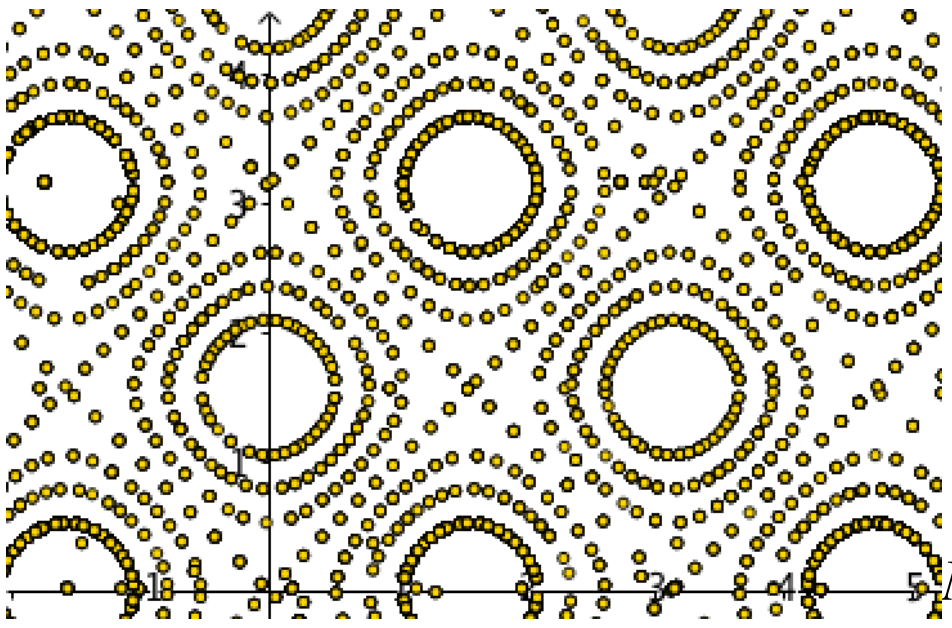

*Preferencialmente deve usar a **Folha Gráfica 2D 2** para estas representações.*

*No caso da função definida em i) só podemos usar o traço sobre* `PMt`.

## *Tarefa 8 - Modelação com Epiciclos*

1. Observe as seguintes três imagens da figura abaixo, estando descrito, ao lado de cada uma das três imagens, o procedimento necessário para as obter no GeoGebra. Foi sempre desenhada uma curva, onde o extremo do domínio de parametrização depende do parâmetro inteiro, `n`, modelado com o comando:

`n=Seletor(1, 12, 1, 1, 200).`

c1 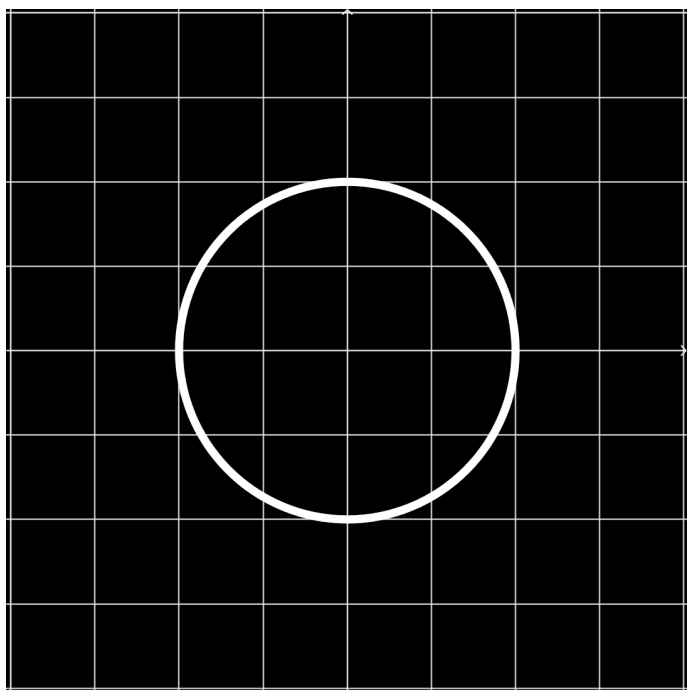 c2 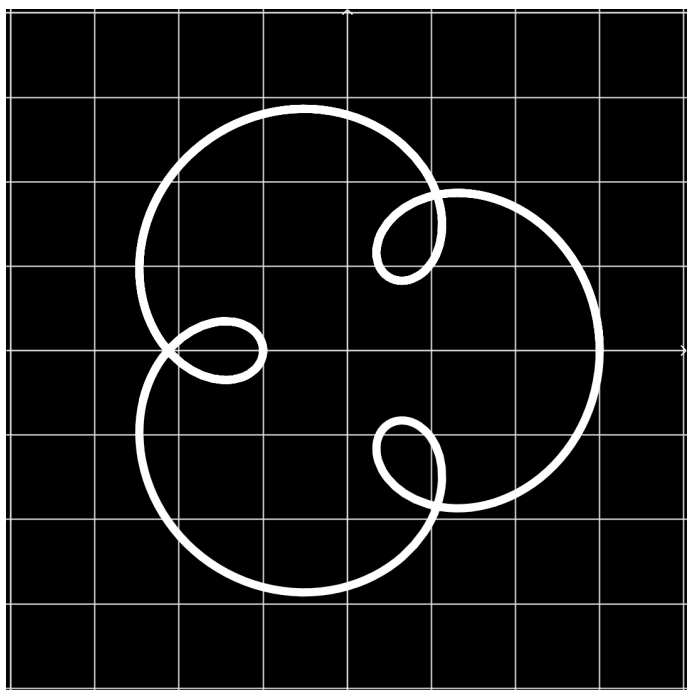 c3 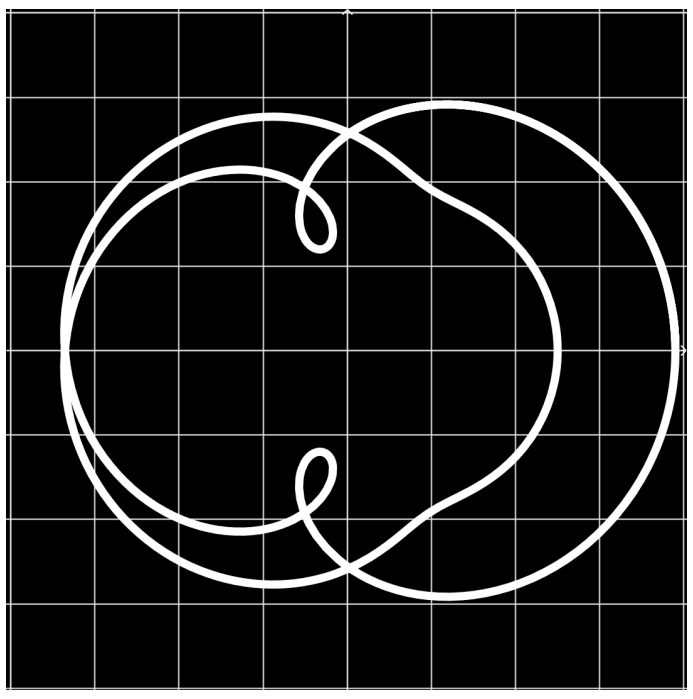

`c1=Curva((2cos(5t/2),2sin(5t/2)),t,0,2nπ)`

`c2=Curva((2cos(2t/5)+cos(8t/5),2sin(2t/5)+sin(8t/5)),t,0,2n π)`

`c3=Curva((2.5 cos(0.2t)+0.7 cos(0.6t)+0.7 cos(t/2),2.5 sin(0.2t)+0.7 sin(0.6t)+ 0.7 sin(t/2)),t,0,2n π)`

   a) Em cada caso, conjeture quais poderão ser os valores do parâmetro `n`.
   b) Obtenha cada uma das imagens no GeoGebra.

2. Disponha o seu ambiente de trabalho no GeoGebra de modo semelhante ao da figura seguinte.

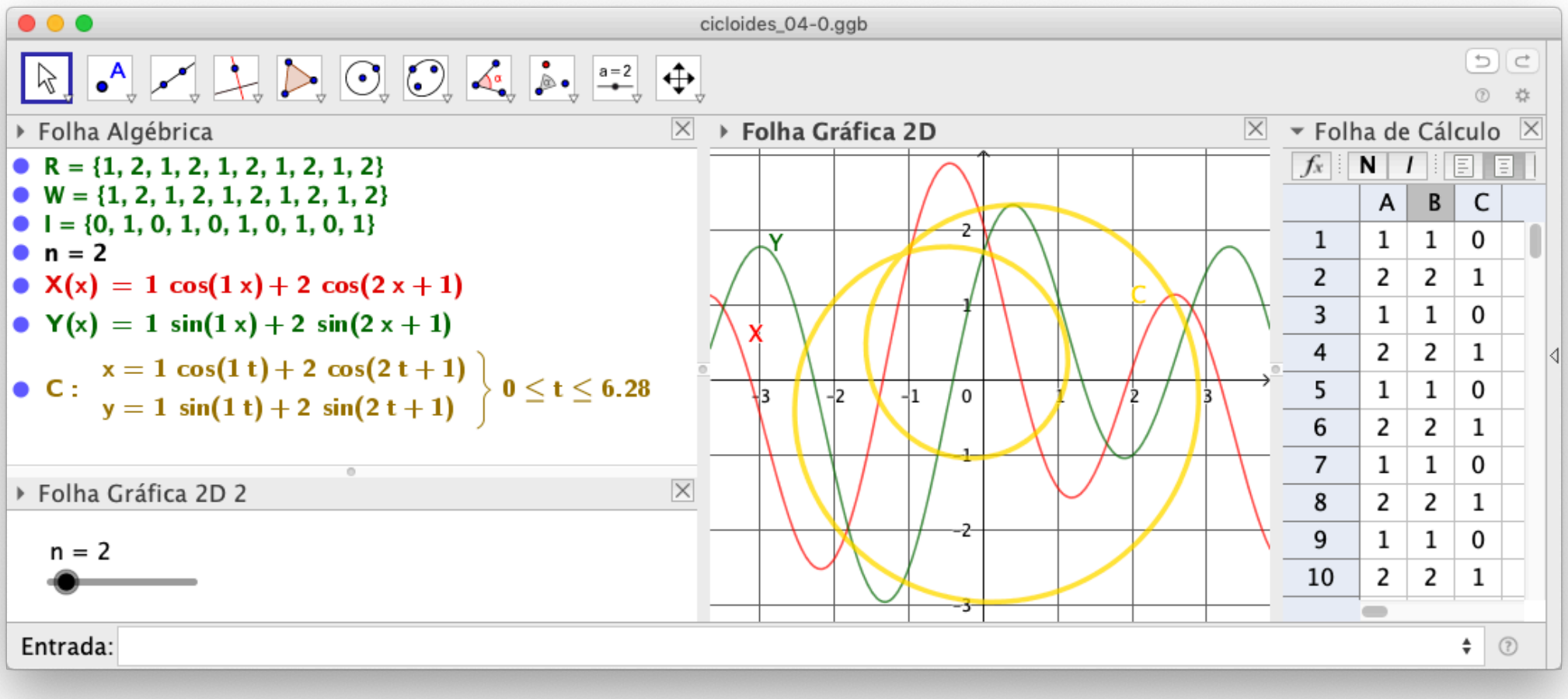

   a) Na **Folha de Cálculo** e insira 10 valores na coluna *A*, *B* e *C*.
   b) Seleccione as células de cada coluna, use a ferramenta ,lista, e defina as listas *R*,*W* e *I* com o conteúdo das células nas colunas A, B e C , respectivamente.
   c) Defina, na **Folha Gráfica 2D 2** `n=Seletor(1,Comprimento(R),1,1,100)` .
   d) Defina, na **Folha Gráfica 2D**:
      i) `X(x)=Soma(Elemento(R,i)*cos(Elemento(W,i)*x+Elemento(I,i)),i,1,n)`
      ii) `Y(x)=Soma(Elemento(R,i)*sin(Elemento(W,i)*x+Elemento(I,i)),i,1,n)`
      iii) `C=Curva(X(t), Y(t), t, 0, 2π )` , na **Folha Gráfica 2D**
   e) Explore a aplicação!

## *Tarefa 9 – Modelação de Airfoils com Números Complexos*

Olive (1997) apresenta um exemplo do uso das transformadas de Joukowski de modo a construir perfis de asas de avião, *airfoils*, a partir de uma circunferência usando o software The Geometer's Sketcpad. A transformada de Joukowski corresponde a função complexa de variável complexa, J, tal que $J(z)=z+\dfrac{1}{z}, z\in\mathbb{C}$.

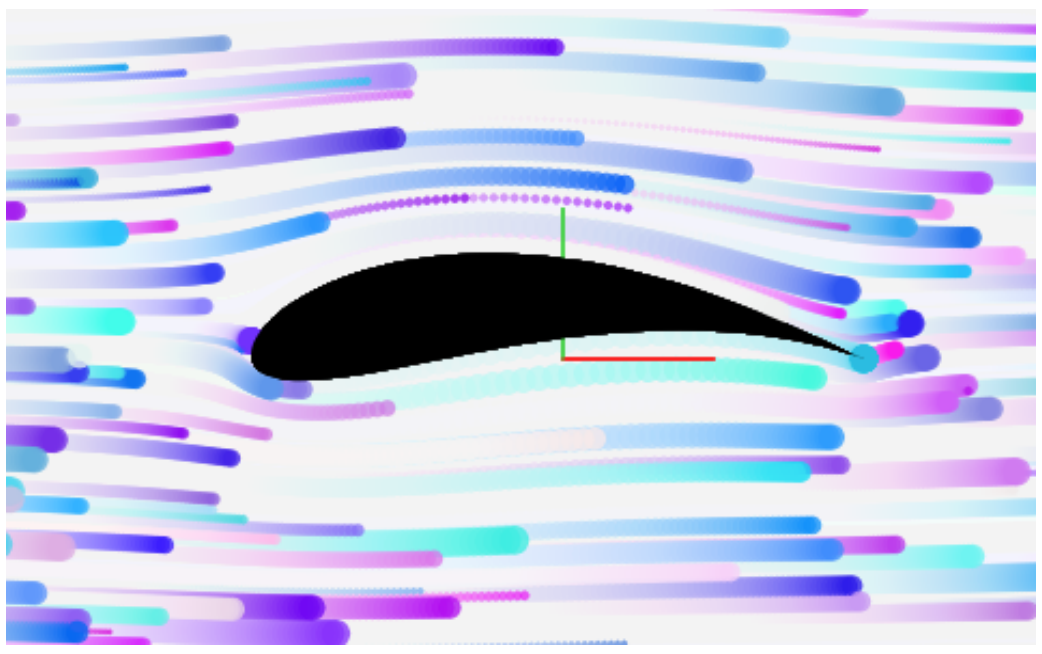

Alerón (Airfoil) (Ponce, 2018 )

1. Abra o GeoGebra, no menu vista ative a ***Folha Gráfica 2 D*** e a ***Folha Gráfica 2 D 2***.
2. Para se assegurar que faz a marcação dos objetos na **Folha Gráfica 2D**, por exemplo, clique com a ferramenta mover, , sobre a **Folha Gráfica 2D**.

3. Com a ferramenta número complexo, , marque um ponto na **Folha Gráfica 2D**. O ponto será designado por $z_1$. Note que: para se referir a este ponto no campo de Entrada terá de usar `z_1`; em alternativa podia ter escrito no campo de *Entrada* `z_1=1+i`, por exemplo.

4. Com a ferramenta número complexo, , marque um ponto na **Folha Gráfica 2D**. O ponto será designado por $z_2$. Note que para se referir a este ponto nos comandos terá de grafar `z_2` .

5. Com a ferramenta circunferência dado o centro e um ponto, , marque a circunferência de centro $z_1$ e que contem $z_2$. A circunferência será designado por `c`. Em alternativa, podia escrever no campo de **Entrada** o comando `Circunferência(z_1,z_2)` pressionando de seguida *Enter*.

6. Com a ferramenta número complexo, , marque um ponto cobre a circunferência . O ponto será designado por $z_3$. Assegure-se que o ponto se move sobre a circunferência , usando a ferramenta .

7. Clique com a ferramenta , mover, sobre a **Folha Gráfica 2D 2**. De seguida, no campo de Entrada introduza o comando z_4=z_3+1/z_3.

8. Com **a Folha Gráfica 2D 2** ativa, use a ferramenta lugar geométrico, , e clique no ponto $z_4$ e de seguida, na **Folha Algébrica**, na representação do ponto $z_3$. Em alternativa, podia introduzir o comando `J=Lugar_Geométrico(z_4, z_3)`, no campo de **Entrada**.

9. Com a ferramenta mover, , modifique a posição dos pontos z1 e z2 na **Folha Gráfica 2D**. Observe as alterações que se verificam na imagem da circunferência pela transformada de Joukowski.

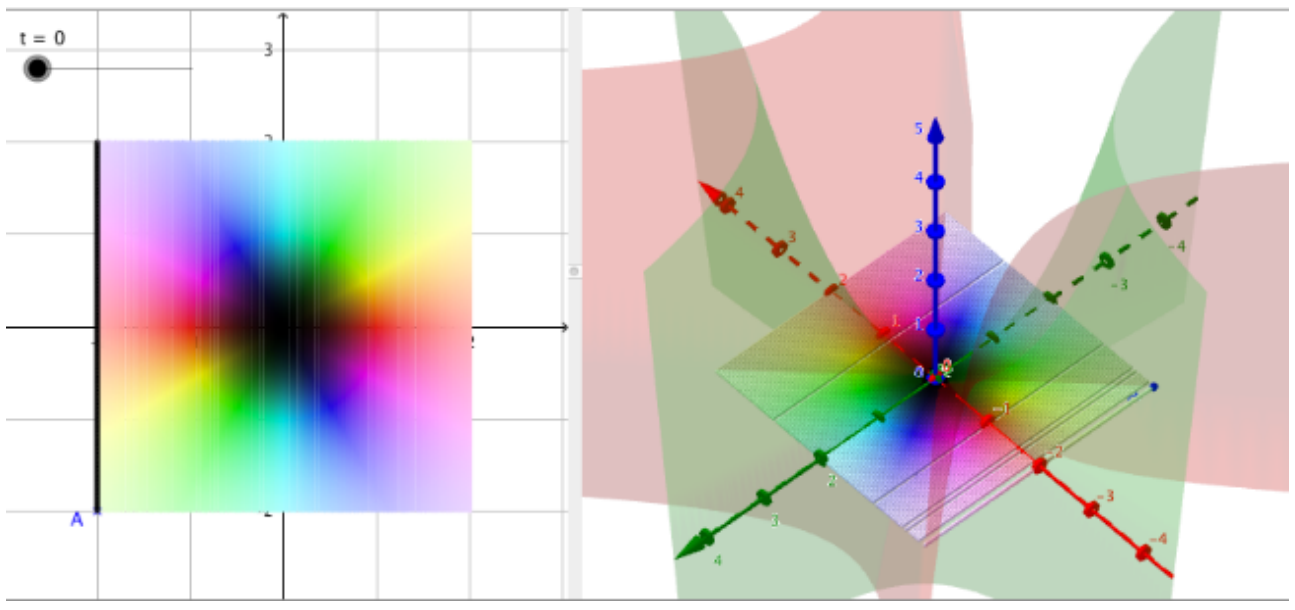

Note que ambas as imagens bidimensionais obtidas são modelos do plano de Argand, representando a **Folha Gráfica 2D** e a **Folha Gráfica 2D 2**, repetivamente, o domínio e o contradomínio da transformada de Joukowski. O uso desta técnica associada aos domínios coloridos permite realizar representações de funções complexas de variáveis complexas, na imagem á esquerda observamos a representação do gráfico de f(z)=z^2, $z\in\mathbb{C}$ e das suas funções componentes. ( Breda & Dos Santos, 2015, 2016; Breda, Dos Santos & Trocado, 2015; Dos Santos & Breda, 2018).

 

## Bibliografia